\documentclass[twocolumn]{autart}    

\usepackage{cite}
\usepackage{amsmath,amssymb,amsfonts}
\usepackage{mathtools}
\usepackage{algorithmic}
\usepackage{graphicx}
\graphicspath{{Figures/}{../Figures/}}
\usepackage{textcomp}
\usepackage{booktabs}
\usepackage{xcolor}
\usepackage{multirow}
\usepackage{cuted}
\usepackage{xstring}
\usepackage{enumerate}
\usepackage{adjustbox}

\DeclareMathOperator*{\E}{\mathbb{E}}

\newcommand{\bt}{\fontseries{b}\selectfont}

\providecommand{\com[2]}{\begin{tt}[#1: #2]\end{tt}}
\definecolor{red}{rgb}{1,0,0}

\definecolor{blu}{rgb}{0,0,1}

\definecolor{gre}{rgb}{0,0.7,0.3}

\newtheorem{definition}{Definition}
\newtheorem{corollary}{Corollary}
\newtheorem{problem}{Problem}

\newcommand{\eem}{~uniformly excited~}
\newcommand{\dm}{~direct module~}
\begin{document}

\begin{frontmatter}

\title{Optimal allocation of excitation and measurement for
identification of dynamic networks \thanksref{footnoteinfo}} 

\thanks[footnoteinfo]{This paper was not presented at any IFAC 
meeting. Corresponding author E.~Mapurunga.}

\author[ASB]{Eduardo Mapurunga}\ead{eduardo\textunderscore mapurunga@ufrgs.br},    
\author[ASB]{Alexandre Sanfelici Bazanella}\ead{bazanella@ufrgs.br}              

\address[ASB]{Data-Driven Control Group \\ Department of Automation and Energy, Universidade Federal do Rio Grande do Sul (DELAE/UFRGS) \\ Porto Alegre-RS, Brazil}             
          

\begin{abstract}                          
	In this paper,  the problem of choosing
	the best allocation of excitations and measurements for the identification
	of a dynamic network is formally stated and analyzed.
	The best choice will be one that achieves the most accurate identification
	with the least costly experiment. Accuracy is assessed by the trace of the
	asymptotic covariance matrix of the parameters estimates, whereas 
	the cost criterion is the number of excitations and measurements.
	Analytical and numerical results are presented for two classes of dynamic
	networks in state space form: branches and cycles. From these results,
	a number of guidelines for the choice emerge, which are based either
	on the topology of the network or on the relative magnitude of the modules being 
	identified.
	An example  is given to illustrate that these guidelines can to some extent
	be applied to networks of more generic topology.  
\end{abstract}

\end{frontmatter}

\section{Introduction}

A crucial aspect of the identification of dynamic networks is the determination 
of which signals must be measured and which nodes must be excited in order
to ensure that the desired modules can be
uniquely identified. This is a fundamental experiment design
problem that has been given a lot
of attention in the recent  literature, and different
facets of it have been attacked.

Verifiable necessary and sufficient conditions 
have been proven for the {\em generic identifiability}
of a whole column (or a whole line) of the
network matrix in \cite{Bazanella2017WhichNodes, Hendrickx2018PartialNodes}.
These conditions concern the measurements of
a fully excited network (or the excitations of a fully measured network),
rely only on the network's topology and can be used to determine the
identifiability of the network as a whole. They have been used in
\cite{Cheng_Shi_VandenHof_2019} to derive a formal method 
for the synthesis of a set of required required excitations in a fully measured network.
An alternative method, that guarantees satisfaction of a slightly different
requirement - that of  {\em global
identifiability} of the network - has been given in \cite{vanwaardeNecessarySufficientTopological2019}. 
Various conditions for generic identifiability have been given, in \cite{Baza2019CDC},
for the more general situation where the network is neither fully measured nor
fully excited, mainly for two particular classes of networks - trees and cycles. 

The problems related to the identification of a particular module in the network
deserve particular attention, and they have received a fair amount of it. 
The identifiability and consistency of the estimates of a single
module have been studied in detail in several publications \cite{materassiIdentificationNetworkComponents2015, 9029448, Dankers_VandenHof_Bombois_Heuberger_2015, WeertsAuto2018,weertsAbstractionsLinearDynamic2020},
resulting in a number of important theoretical concepts and in identification methods \cite{VanDenHof2013ComplexPem, Dankers2016PredInputSel, geversPracticalMethodConsistent2018}
involving different choices of measurements and excitations. 
Correlation among the noise sources affecting the different nodes of the network
poses additional difficulties to achieve identifiability and consistency, and this issue
has also been dealt with in the literature, with identifiability analysis provided in
\cite{geversIdentifiabilityDynamicalNetworks2019} and identification methods provided in \cite{weertsPredictionErrorIdentification2018, everittEmpiricalBayesApproach2018}.

Whether the whole network is to be identified, or only a small part of it - possibly
a single module - is of interest, usually there are multiple
choices for the excitations and measurements that guarantee identifiability. 
Once identifiability is achieved,
accuracy of the identification from noisy signals becomes naturally the next matter of concern in the design of an experiment for the identification of a network.
Specifically, one would like to choose the nodes to be excited  and those to
measure such that the accuracy of the identification is optimized, always
employing the minimum number of measurements and excitations possible. 
This problem is what this paper is about. 
We first provide a theoretical framework in which to study the problem, then present 
theoretical analysis and numerical studies for two classes of networks: 
branches and cycles where the edges consist of single delay transfer functions. 

This paper is organized as follows.
The dynamic network system setup, which is standard in the literature, is briefly presented in Section \ref{sec:Dynetsetting},
along with the particular classes of networks - cycles and branches - that are studied in depth in this paper. 
Then, in Section \ref{sec:Problem}, we formally state the problem, which requires a number of new mathematical
definitions and a review of some known but not standard concepts and results on generic identifiability. 
The main results of the paper appear in Sections  \ref{sec:Loops} and \ref{sec:Branches}, 
where variance analysis is carried out for cycles and for branches, respectively.  
An example combining these two structures is presented in Section \ref{sec:anexample} to verify 
whether the results obtained for cycles and branches can be applied to more general topologies.
Concluding remarks are given in Section \ref{sec:Conclusion}.

\section{Dynamic Networks Setting}
\label{sec:Dynetsetting}

We consider dynamic networks composed of $n$ nodes which represent scalar internal signals $\left\lbrace w_j \right\rbrace$ for $j \in \mathcal{W} \triangleq \left\lbrace 1, 2, \dots, n\right\rbrace$. 
The nodes may have the influence of external signals $\left\lbrace r_j(t) \right\rbrace$ and
the  internal signals are accessible through a noisy measurement. 
Such dynamic networks can be described by the following network equation:
\begin{subequations}
	\begin{align}
	w(t) &= G^0(q)w(t) + Br(t), \label{eq:S_internal}\\
	y(t) &= Cw(t) + e(t), \label{eq:S_external}
	\end{align}
\end{subequations}
where 
$y(t) \in \mathbb{R}^p$ is the network output,
$w(t) \in \mathbb{R}^n$, $r(t) \in \mathbb{R}^m $, $e(t) \in \mathbb{R}^p$ are the network's internal and external signals and the corrupting measurement noise, respectively.
The matrix $G^0(q)$ is referred as network matrix, whose entries are causal discrete-time rational transfer functions in the operator $q$, the  forward shift operator, i.e. $qr(t) = r(t+1)$.
We refer to (\ref{eq:S_internal}) as the network equation and to  (\ref{eq:S_external}) as the measurement equation.
The matrices $B \in \mathbb{Z}_2^{n \times m}$ and $C \in \mathbb{Z}_2^{p \times n}$ are binary selection matrices with full row rank and full column rank, respectively. They play the role of selecting which external signals are applied to the nodes and which measurements are taken from the network. Associated with these matrices, we will additionally denote the set of excited nodes by $\mathcal{B}$ and measured nodes by $\mathcal{C}$ \cite{Hendrickx2018PartialNodes,Bazanella2017WhichNodes,Baza2019CDC}.

In this paper we pose the problem of excitation and measurement allocation for general networks in the form
(\ref{eq:S_internal})(\ref{eq:S_external}), then concentrate our analysis on two particular classes of networks:
cycles and branches.
A cycle - or loop - is a connected network such that there's exactly one path from any particular node to itself.
A branch, on the other hand, is a connected network such that there exists exactly one path from
node $i$ to node $j$ if and only if $j > i$. 
Conditions for generic identifiability of these two classes of networks are known, and these conditions can also be applied for these structures when they form subgraphs
of a larger network \cite{Baza2019CDC}, provided some additional conditions on the graph are satisfied.

In order to arrive at clear, intuitive conditions, the analysis in this paper is developed for  
cycles and branches in which the network matrix has the form $G^0(q) = A^0q^{-1}$
and the entries of $A^0$ are real numbers. It is easy to see that the network equation thus obtained
constitutes a state space equation, and accordingly a variety of results in network theory have
been presented for state-space equations \cite{vanwaardeIdentifiabilityUndirectedDynamical2018}. 

We shall thus study cycles whose network matrices have the following form:
\begin{equation}
G^0 (q)= q^{-1} \begin{bmatrix}
0 & 0 & 0 & \cdots & a^0_{1n}\\
a^0_{21} & 0 & 0 & \cdots & 0\\
0 & a^0_{32} & \ddots & \cdots & \vdots\\
\vdots &\cdots & \ddots & \ddots & \vdots\\
0 & \cdots & \cdots & a^0_{n,n-1} & 0\\
\end{bmatrix}\label{eq:GLoop},
\end{equation}
where $a^0_{ji}$ are real numbers and we have numbered the nodes sequentially
(that is, $a_{ji}\neq  0$ if and only if $j=i+1$ or $j=1$ and $i=n$) without loss of generality. 
The network matrix of a branch, in its turn, has the following form:
\begin{equation}
G^0 (q)= q^{-1} \begin{bmatrix}
0 & 0 & 0 & \cdots & 0\\
a^0_{21} & 0 & 0 & \cdots & 0\\
0 & a^0_{32} & \ddots & \cdots & \vdots\\
\vdots &\cdots & \ddots & \ddots & \vdots\\
0 & \cdots & \cdots & a^0_{n,n-1} & 0\\
\end{bmatrix}\label{eq:GCasc},
\end{equation}
where again we have numbered the nodes sequentially (that is, $a_{ji}\neq  0$ if and only if $j=i+1$)
without loss of generality. In all the discussion to follow it is assumed this
sequential numbering of the nodes for both classes of networks.

In analysing these two classes of networks with simple transfer functions we will be able to provide 
closed-form solutions for the problem at hand and derive the concepts behind these solutions. 


\section{Problem statement}
\label{sec:Problem}

In this Section, we first review identifiability conditions for branches and cycles.
Then the problem is formally stated and the experimental conditions and some technical assumptions
are presented. 

Identifiability of a network depends on the graph associated with the network matrix $G^0(q)$ and on the two selection matrices $B$ and $C$. For a given network
matrix $G^0(q)$, different combinations of selection matrices can provide identifiability, and one would like to use the smallest number of excitations and/or measurements; this motivates the following definition.

\begin{definition}
	Let  $G^0(q)$ be a network matrix, for which different selection matrices are considered.
	A pair of selection matrices $B$ and $C$, with its corresponding node sets $\mathcal{B}$ and $\mathcal{C}$,
	is called an {\em excitation and measurement pattern} - EMP, for short. An EMP is said to be a {\em valid}
	EMP if it is such that the network (\ref{eq:S_internal})-(\ref{eq:S_external}) is generically identifiable.
	Define $\nu = \vert \mathcal{B} \vert +
	\vert \mathcal{C} \vert$\footnote{$\vert \cdot \vert $ standing for the cardinality of a set.} as the {\em cardinality} of an EMP.
	A given EMP is said to be a {\em minimal EMP} if it is valid and there is no valid EMP with smaller cardinality. 
	\label{def:EMP}
\end{definition}

The following result gives a necessary condition for generic identifiability of a general network - stated otherwise,
for the validity of an EMP.
\begin{thm}\cite{Baza2019CDC}
	All transfer functions in a network are generically identifiable only if ${\mathcal B} \cup {\mathcal C} = {\mathcal W}$. 
	\label{teo:ExciteOrMeasure}
\end{thm}
In other words, any valid EMP must contain all nodes. For branches one does not need much more than this to guarantee
identifiability, as stated formally in the following Corollary - which is an obvious consequence of Theorem V.1 in \cite{Baza2019CDC}.

\begin{corollary} 
	A branch is generically identifiable if and only if the following conditions are satisfied:
	\begin{enumerate}
		\item $1 \in {\mathcal B}$ (the source node is excited);
		\item $n \in {\mathcal C}$ (the sink node is measured);
		\item $\mathcal{B}\, \cup\, \mathcal{C} = {\mathcal W}$ (all nodes are either measured or excited).
	\end{enumerate}
\end{corollary} 

For dynamic networks of the cycle type, we have another condition regarding the minimum number of excitation/measurement setting.

\begin{thm}\cite{Baza2019CDC}
	A  cycle is generically
	identifiable if $\mathcal{B} \cup \mathcal{C} = \mathcal{W}$ and $\mathcal{B} \cap \mathcal{C} \neq \emptyset$.
	\label{teo:IdentLoop}
\end{thm}

This result gives a sufficient condition for generic identifiability of a whole cyclic network, on top of the necessary condition in Theorem
\ref{teo:ExciteOrMeasure}:  at least one node must be both excited and measured. This is not a necessary
condition, as it can be relaxed  for some network classes with some particular excitation/measurement patterns, as in the following theorem.

\begin{thm}\cite{Baza2019CDC}
	For cyclic dynamic network with $n > 3$ nodes and $n$ even, all transfer functions can be generically identified if its nodes are alternately measured and excited.
	\label{teo:CyclicEven} 
\end{thm}  

As a direct consequence of these known results and the above definitions, we have the following. 

\begin{corollary} \label{cor:summary}
	For a branch, an EMP is minimal if and only if it satisfies the four conditions:
	\begin{enumerate}
		\item $1\in {\mathcal B}$
		\item $n\in {\mathcal C}$
		\item ${\mathcal B}\cup {\mathcal C}={\mathcal W}$
		\item ${\mathcal B}\cap {\mathcal C}=\emptyset$.
	\end{enumerate}
	As a consequence, there are $2^{n-2}$  minimal EMPs in a branch.
	
	For a cycle, an EMP is minimal if it satisfies the following:
	\begin{itemize}
		\item for $n$ even and larger than $2$, even nodes are excited and odd nodes are measured, or vice-versa;
		\item for other $n$, $\mathcal{B} \cup \mathcal{C} = \mathcal{W}$ and $\mathcal{B} \cap \mathcal{C} =\{i\} $
		for some $i \in \mathcal{W}$.
	\end{itemize}
	Therefore, in cycles with an even number of nodes there are only $2$ minimal EMPs, whereas for odd number of nodes there are
	$n\times 2^{n-1}$ minimal EMPs.
\end{corollary}

{\bf Remark} - Please note that the conditions  for cycles in Corollary \ref{cor:summary}
are sufficient but not necessary. So, when an EMP satisfies one of these
(sufficient) conditions it is not necessarily a minimal EMP, as per Definition \ref{def:EMP}.
Yet it is, under the current state of the theory, the EMP with smallest cardinality that we can guarantee to be valid.
For this reason we will say that such EMPs are minimal, even though we can not establish for sure that they are. 

In this paper we are interested in an experiment design question, which is deciding which minimal EMP is best, that is, the one that provides the most accurate identification. 
We will assess accuracy through the  Cramer-Rao lower bound, which is given by the asymptotic parameter covariance
matrix obtained by  prediction error identification. This matrix can be calculated as
$P = [\E\psi(t)\Lambda^{-1}\psi(t)^T]^{-1}$ \cite{ljung1998system}, where $\psi(t)$ is the gradient of the
optimal one-step ahead predictor and $\Lambda$ is the noise's covariance matrix. The overall accuracy
is assessed by the trace of the matrix $P$, a criterion
known in experiment design literature as A-optimality \cite{pukelsheimOptimalDesignExperiments2006}.

Regarding the experimental setup, we will consider that the following assumptions hold within the prediction error framework:
\begin{enumerate}[(a)]
	\item the dynamic network is stable; 
	\item the external signals $\lbrace r_j(t)\rbrace$ are mutually independent stationary white noise processes with zero mean and variance $\sigma_j^2$; they are uncorrelated with all noise processes $\left\lbrace e_j\right\rbrace$; \label{item:input-white-noise}
	\item the corrupting noise sequences $\left\lbrace e_j\right\rbrace $ are independent stationary Gaussian white noise processes with zero mean and variance $\lambda_j$.  \label{item:noise-lambda}
\end{enumerate}
We can now formally state the problem studied in this paper.

\begin{problem}
	Given a dynamic network with cycle or branch structure and satisfying the technical conditions (a), (b), (c) above,
	determine what is the minimal EMP 
	that results in the smallest value for the trace of the asymptotic covariance matrix $P$. 
\end{problem}


Before presenting the analysis for these two classes of dynamic networks, we will define additional nomenclature 
that will be used throughout this work.
\begin{definition}
	A module $a^0_{ji}q^{-1}$ is said to be a \emph{direct module} of an EMP if
	$i \in \mathcal{B}$ and $j \in \mathcal{C}$.
	\label{def:directmodule}
\end{definition}
The purpose of this definition will become clear later on, where this kind of module will play a key role in determining which minimal EMP is the most accurate. 
The next definition is related to the experimental setting for which the network is submitted. 
\begin{definition}
	A dynamic network is said to be \emph{uniformly excited} if $\lambda_j = \lambda$ and  $\sigma_i^2 = \sigma^2$ for all $j \in \mathcal{C}$ and $i \in \mathcal{B}$.
	It is  a \emph{fully symmetric} network when, in addition to being uniformly excited,
	all modules of the network have the same magnitude.
	\label{def:UE} 
\end{definition}

\section{Cycles}
\label{sec:Loops}

In this section we present the analysis for dynamic networks with cycle structure, which are characterized
by the network matrix in (\ref{eq:GLoop}).
We start with the simplest case, with only two nodes, and build on it by adding more nodes to the
analysis.

\subsection{Two-nodes}
\label{secsec:2nLoops}

In two-node cycle there are four minimal EMPs, namely: 
\begin{enumerate}[(I)]
	\item
	$ \mathcal{B} = \left\lbrace 1, 2 \right\rbrace; \mathcal{C} = \left\lbrace 1\right\rbrace$, \label{item:2cEMP1}
	\item
	$ \mathcal{B} = \left\lbrace 1, 2 \right\rbrace; \mathcal{C} = \left\lbrace 2\right\rbrace$, \label{item:2cEMP2}
	\item
	$ \mathcal{B} = \left\lbrace 1\right\rbrace; \mathcal{C} = \left\lbrace 1, 2\right\rbrace$, \label{item:2cEMP3}
	\item
	$ \mathcal{B} = \left\lbrace 2 \right\rbrace; \mathcal{C} = \left\lbrace 1, 2\right\rbrace$. \label{item:2cEMP4}
\end{enumerate}

In order to calculate the asymptotic covariance matrix of the parameter estimates one needs the gradient of the optimal predictor. 
Recall that the noise processes obey assumption (\ref{item:noise-lambda}).
The optimal predictors, in this case, are as follows:
\begin{align}
\hat{y}_1(t|t-1) &= ({r_1(t) + a_{12}r_2(t-1)})/({1 - a_{12}a_{21}q^{-2}}), \nonumber\\
\hat{y}_2(t|t-1) &= ({r_2(t) + a_{21}r_1(t-1)})/({1 - a_{12}a_{21}q^{-2}}). \nonumber
\end{align}
Hence, the gradient of the optimal predictors in EMP (I) will be:
\begin{align}
\psi^{\text{I}}(t)^T = \frac{1}{\Delta}
\begin{bmatrix}
r_2(t-1) + a_{21}r_1(t-2)  \\
a_{12}\left(r_1(t-2) + a_{12}r_2(t-3)\right) 
\end{bmatrix}, \nonumber
\end{align}
where 
$
\Delta = 1 - 2a_{12}a_{21}q^{-2} + a_{12}^2a_{21}^2q^{-4}
$. 
Let $\hat{a}_{12}^{\text{I}}$ denote an estimate of the true parameter $a_{12}^0$, with the superscript indicating the EMP. Since the network's inputs obey assumption (\ref{item:input-white-noise}), one can obtain the variance of the parameter estimates by computing $\E\psi^{\text{I}}(t)\psi^{\text{I}}(t)^T/\lambda_1$, yielding:
\begin{align}
\text{var}(\hat{a}_{12}^{\text{I}}) &= {\lambda_1\gamma_0}(\sigma_1^2 + \sigma_2^2[a^0_{12}]^2)/{d_1}, \label{eq:Ex2x2vb121} \\
\text{var}(\hat{a}_{21}^{\text{I}}) &= {\lambda_1\gamma_0}{( [a^0_{21}]^2\sigma_1^2 + \sigma_2^2)/([a^0_{12}]^2d_1}), \label{eq:Ex2x2bb211}
\end{align}
where 
$d_1 = \left( \gamma_0^2M_{\sigma_1^2, \sigma_2^2}N_{\sigma_1^2, \sigma_2^2} - \left(\gamma_0^2\sigma_1^2a^0_{21} + \gamma_2\sigma_2^2a^0_{12}\right)^2 \right)$, 
$M_{\sigma_1^2, \sigma_2^2} \triangleq \left( \sigma_1^2[a^0_{21}]^2 + \sigma_2^2\right) $,
$N_{\sigma_1^2, \sigma_2^2} \triangleq \left( \sigma_1^2 + \sigma_2^2[a^0_{12}]^2\right)  $,
and 
\begin{align*}
\gamma_0 &= \frac{[a_{12}^0 a_{21}^0]^{4} + 8 [a_{12}^0 a_{21}^0]^{2} + 1}{	[a_{12}^0 a_{21}^0]^{8} + 4 [a_{12}^0 a_{21}^0]^{6} - 6 [a_{12}^0 a_{21}^0]^{4} + 4 [a_{12}^0 a_{21}^0]^{2} + 1},\\
\gamma_2 &= \frac{2 a_{12}^0 a_{21}^0 \left( [a_{12}^0 a_{21}^0]^{2} + 1\right) }{[a_{12}^0 a_{21}^0]^{8} + 4 [a_{12}^0 a_{21}^0]^{6} - 6 [a_{12}^0 a_{21}^0]^{4} + 4 [a_{12}^0 a_{21}^0]^{2} + 1}.
\end{align*} 
Now, in EMP (II) we have the following result:
\begin{align}
\psi^{\text{II}}(t) = \frac{1}{\Delta}\begin{bmatrix}
a_{21}\left(a_{21}r_1(t-3) + r_2(t-2) \right) \\
r_1(t-1) + a_{12}r_2(t-2)
\end{bmatrix}. \nonumber
\end{align}
Similarly, the variances of the parameters' estimates are:
\begin{align}
\text{var}(\hat{a}_{12}^{\text{II}}) &= {\lambda_2\gamma_0}(\sigma_1^2 + \sigma_2^2[a^0_{12}]^2)/({[a^0_{21}]^2 d_2}), \label{eq:Ex2x2bb122}\\
\text{var}(\hat{a}_{21}^{\text{II}}) &= {\lambda_2\gamma_0}([a^0_{21}]^2\sigma_1^2 + \sigma_2^2)/{d_2},\label{eq:Ex2x2bb212}
\end{align}
where
$d_2 = \left(\gamma_0^2M_{\sigma_1^2, \sigma_2^2}N_{\sigma_1^2, \sigma_2^2} - \left(\gamma_0^2\sigma_2^2a^0_{21} + \gamma_2\sigma_1^2a^0_{12}\right)^2 \right)$. 
When EMP (III) is considered, the gradient of the optimal predictor has the following form:
\begin{align}
\psi^{\text{III}}(t) = \frac{1}{\Delta}\begin{bmatrix}
{a_{21}r_1(t-2)} & {a_{21}^2r_1(t-3)} \\
{a_{12}r_1(t-2)} & {r_1(t-1)}
\end{bmatrix}. \nonumber
\end{align}
Thus, the variances of the parameters' estimates are: 
\begin{align}
\text{var}(\hat{a}_{12}^{\text{III}}) &= {\lambda_1\lambda_2\gamma_0}(\lambda_1 + \lambda_2[a^0_{12}]^2)/({[a^0_{21}]^2d_3}), \label{eq:Ex2x2bb123}\\
\text{var}(\hat{a}_{21}^{\text{III}}) &= {\lambda_1\lambda_2\gamma_0}([a^0_{21}]^2\lambda_1 + \lambda_2 )/{d_3},
\end{align}
with $d_3 =  \sigma_1^2\left(\gamma_0^2M_{\lambda_1, \lambda_2}N_{\lambda_1, \lambda_2} - \left(\gamma_0^2\lambda_2a^0_{12} + \gamma_2\lambda_1a^0_{21}\right)^2 \right)$.
Finally, the last EMP (IV) leads to the following predictor gradient:
\begin{align}
\psi^{\text{IV}}(t) = \frac{1}{\Delta}\begin{bmatrix}
{r_2(t-1)} & {a_{21}r_2(t-2)} \\
{a_{12}^2r_2(t-3)} & {a_{12}r_2(t-2)}
\end{bmatrix}, \nonumber
\end{align}
which leads to the following variances of the parameters' estimates:
\begin{align}
\text{var}(\hat{a}_{12}^{\text{IV}}) &= {\lambda_1\lambda_2\gamma_0}(\lambda_1 + \lambda_2[a^0_{12}]^2)/{d_4}, \\
\text{var}(\hat{a}_{21}^{\text{IV}}) &= {\lambda_1\lambda_2\gamma_0}([a^0_{21}]^2\lambda_1 + \lambda_2)/({[a^0_{12}]^2d_4}), \label{eq:Ex2x2vb214}
\end{align}
with $d_4 =  \sigma_2^2\left(\gamma_0^2M_{\lambda_1, \lambda_2}N_{\lambda_1, \lambda_2} - \left(\gamma_0^2\lambda_2a^0_{12} + \gamma_2\lambda_1a^0_{21}\right)^2 \right)$.

A number of conclusions can be drawn from these equations. First consider a fully symmetrical network, that is, one in which $a^0_{12} = a^{0}_{21}$and that is
uniformly excited (see Definition \ref{def:UE}).

Let $P^j$ denote the covariance matrix of EMP $j$; the analytical expressions
for these $P^j$'s are given in Appendix \ref{app1}, from which the following facts
are observed: 
\begin{enumerate}
	\item $P^{\text{I}} = P^{\text{IV}}$ and $P^{\text{II}} = P^{\text{III}}$; \label{item:2nodesymmetric-CovsEqual}
	\item $tr(P^{\text{I}}) = tr(P^{\text{II}})= tr(P^{\text{III}})= tr(P^{\text{IV}})$; \label{item:2nodesymmetric-trEqual}
	\item $var(\hat{a}^{\text{I}}_{12}) < var(\hat{a}^{\text{II}}_{12})$, $var(\hat{a}^{\text{I}}_{21}) > var(\hat{a}^{\text{II}}_{21})$, $var(\hat{a}^{\text{IV}}_{12}) < var(\hat{a}^{\text{III}}_{12})$ and $var(\hat{a}^{\text{IV}}_{21}) > var(\hat{a}^{\text{III}}_{21})$.
\end{enumerate}
The overall precision, given by the trace of the covariance matrix,
is the same in all four EMPs (condition 
(\ref{item:2nodesymmetric-trEqual})), so from this point of view
it is irrelevant whether one has two measures or two excitations.
Furthermore, EMPs (I) and (IV) result in better accuracy for $\hat{a}_{12}$,
while EMPs (II) and (III) yield a more accurate estimate for $\hat{a}_{21}$; 
each module is estimated with better precision when it is a {\em direct module} 
(see Definition \ref{def:directmodule}).

Now, consider the situation where the network is uniformly excited but with 
arbitrary values of the real parameters. In this scenario, all covariance matrices are inversely proportional to the signal-to-noise ratio (SNR) $\sigma^2/\lambda$.
Then, the expressions for the covariance matrices allow to state the following facts:
\begin{enumerate}
	\item $P^{\text{I}} = P^{\text{IV}}$ and $P^{\text{II}} = P^{\text{III}}$; 
	\item if $|a^0_{21}| > |a^0_{12}|$ then $tr(P^{\text{I}}) > tr(P^{\text{II}})$;
	\item if $|a^0_{21}| < |a^0_{12}|$ then $tr(P^{\text{I}}) < tr(P^{\text{II}})$. 
\end{enumerate}

As in the fully symmetric case, EMPs (I) and (IV) share the same covariance matrix and similarly EMPs (II) and (III). 
Moreover, the largest module has a direct influence in the accuracy of the parameters estimates in the different EMPs, which can be seen by comparing the sum of (\ref{eq:Ex2x2vb121})-(\ref{eq:Ex2x2bb211}) and (\ref{eq:Ex2x2bb122})-(\ref{eq:Ex2x2bb212}). Specifically, smaller trace of the covariance matrix is obtained by the EMPs in which
the largest parameter value is a direct module.

In order to illustrate these results numerically for the \eem~scenario, we consider three representative cases: 1) $a^0_{21} = a^0_{12} = 0.5$, 2) $a^0_{21} = 1$; $a^0_{12} = 0.5$, 3) $a^0_{21} = 0.5$; $a^0_{12 }= 1$.
The theoretical variances obtained from the analytical expressions (\ref{eq:Ex2x2vb121}) -- (\ref{eq:Ex2x2vb214}) are given in Table \ref{tab:2nodeExperiments} 
(multiplied by the SNR) for all minimal
EMPs. It is observed in this Table that the accuracy gains obtained in choosing the best EMPs are very significant.
Indeed, the trace of the covariance matrix for the best EMPs (II and III for case 2, I and IV for case 3), is about five times
smaller than for the other EMPs ($1.06$ versus $4.86$) and the gains obtained in the precision of individual modules get
even larger
($0.65$  versus $4.20$ in the estimate  of $a_{21}$ in Case 2, for example).

\begin{table}[!h]
	\centering
	\caption{Theoretical Variances $\times$ SNR for three cases, all EMPs}
	\label{tab:2nodeExperiments}
	\begin{adjustbox}{width=\columnwidth}
	\begin{tabular}{c| c c c c c c c c}
		\hline
		Case & $\hat{a}^{\text{I}}_{12}$& $\hat{a}^{\text{I}}_{21}$& $\hat{a}^{\text{II}}_{12}$& $\hat{a}^{\text{II}}_{21}$& $\hat{a}^{\text{III}}_{12}$& $\hat{a}^{\text{III}}_{21}$& $\hat{a}^{\text{IV}}_{12}$& $\hat{a}^{\text{IV}}_{21}$\\\hline
		1        & 0.92 & 3.64 & 3.64 & 0.92 & 3.64 & 0.92 & 0.92 & 3.64\\
		2        & 0.66 & 4.20 & 0.41 & 0.65 & 0.41 & 0.65 & 0.66 & 4.20\\
		3        & 0.65 & 0.41 & 4.20 & 0.66 & 4.20 & 0.66 & 0.65 & 0.41 \\
		\hline
	\end{tabular}
	\end{adjustbox}
\end{table}

To summarize the results obtained so far, we have shown
the role played by two factors in the accuracy of the identification:
the existence of direct modules and the relative magnitudes of the modules.
We demonstrated analytically the 
following principles: 
\begin{itemize}
	\item all factors being equal (fully symmetrical case), direct modules are estimated more accurately;
	\item  regarding the modules' magnitudes, all other factors being equal (\eem),
	better overall accuracy is obtained when the largest module
	is a direct module.
\end{itemize}
When and how these principles generalize to other network topologies 
is explored along the paper. We have also illustrated numerically that the gains in precision that can be obtained by choosing the best EMP are quite significant.


\subsection{Three-nodes Cycles}
\label{subsec:3nLoop}

Let us extend our analysis to cyclic dynamic networks with 3 nodes.
There are more choices for the excitation/measurement 
patterns in this case: at least one node must be both excited and measured,
and for each node that is
both excited and measured four different patterns exist; this produces a total of 12 minimal EMPs, which are listed in Table \ref{tab:3nodeScenarios}. 
\begin{table}[!h]
	\caption{The 12 EMPs for a 3-node cycles.}
	\label{tab:3nodeScenarios}
	\centering
	\begin{adjustbox}{width=\columnwidth}
	\begin{tabular}{c c c c}
		\hline
		Nº & EMP & Nº & EMP \\\hline
		I & $\mathcal{B} = \left\lbrace 1, 2, 3 \right\rbrace; \mathcal{C} = \left\lbrace 1\right\rbrace $ & II & $\mathcal{B} = \left\lbrace 1, 2, 3 \right\rbrace; \mathcal{C} = \left\lbrace 2\right\rbrace $ \\
		III & $\mathcal{B} = \left\lbrace 1, 2, 3 \right\rbrace; \mathcal{C} = \left\lbrace 3\right\rbrace $ & IV & $\mathcal{B} = \left\lbrace 1 \right\rbrace; \mathcal{C} = \left\lbrace 1, 2, 3\right\rbrace $ \\
		V & $\mathcal{B} = \left\lbrace 2 \right\rbrace; \mathcal{C} = \left\lbrace 1, 2, 3\right\rbrace $ & VI & $\mathcal{B} = \left\lbrace 3 \right\rbrace; \mathcal{C} = \left\lbrace 1, 2, 3\right\rbrace $ \\
		VII & $\mathcal{B} = \left\lbrace 1, 2 \right\rbrace; \mathcal{C} = \left\lbrace 1, 3\right\rbrace $ & VIII & $\mathcal{B} = \left\lbrace 1, 3 \right\rbrace; \mathcal{C} = \left\lbrace 1, 2\right\rbrace $ \\
		IX & $\mathcal{B} = \left\lbrace 2, 3 \right\rbrace; \mathcal{C} = \left\lbrace 1, 2\right\rbrace $ & X & $\mathcal{B} = \left\lbrace 1, 2 \right\rbrace; \mathcal{C} = \left\lbrace 2, 3\right\rbrace $ \\
		XI & $\mathcal{B} = \left\lbrace 1, 3 \right\rbrace; \mathcal{C} = \left\lbrace 2, 3\right\rbrace $ & XII & $\mathcal{B} = \left\lbrace 2, 3 \right\rbrace; \mathcal{C} = \left\lbrace 1, 3\right\rbrace $ \\\hline
	\end{tabular}
	\end{adjustbox}
\end{table}

This large number of minimal EMPs prevents a meaningful analytical comparison among them in the spirit
of the one performed for the 2-node network. 
So, in order to gain insight in the choice of the best EMP, we have randomly generated $1,000$ case studies
with different parameter values, with each parameter sampled from a uniform distribution with support in $[-1~,~+1]$.
We considered an \eem~scenario, and for each case study the theoretical 
covariance matrix was computed and 
the best minimal EMP was selected as the one which resulted in the
smallest trace of the covariance matrix.
The number of cases in which each EMP was selected in this experiment are shown in Table \ref{tab:3nodeHistogram}.

\begin{table}[!h]
	\caption{Number of times that each EMP was selected as the best from 1,000 randomly generated systems.}
	\label{tab:3nodeHistogram}
	\centering
	\begin{tabular}{c c c c c c c}
		\hline
		EMP    & VII & VIII & IX & X & XI & XII \\\hline
		Frequency &236&117&201&105 &238 &103 \\\hline
	\end{tabular}
\end{table}
EMPs I to VI do not appear in this Table because their count was zero - they were never selected as the best EMP.  What the EMPs I to VI have in common is that in all of them the measurements and excitations are oddly distributed among
the nodes - that is, the cardinality of the sets $\mathcal B$ and $\mathcal C$ are different. 
On the other hand, $|B| = |C|$ for all EMP's that appear as ``winners" in Table
\ref{tab:3nodeHistogram} - EMPs VII to XII.

From the results with 2-node cyclic networks, we conjecture that EMPs for which the largest parameters are direct modules
should result in the smaller trace of covariance matrix. 
In order to test this hypothesis we consider 7 different case studies with different parameter values, described in Table \ref{tab:3nodeParameters}.
The trace of the covariance matrix for EMPs VII-XII in the different cases are displayed in Table \ref{tab:3nodeResults} with $\sigma^2_j = 1$ and $\lambda_j = 0.01$, $\forall j \in \mathcal{W}$, boldface letters indicating the best EMP for each case.
\begin{table}[!h]
	\caption{Parameter values for the numerical experiments.}
	\label{tab:3nodeParameters}
	\centering
	\begin{tabular}{c c c c c c c c}
		\hline
		Exp. & $a^0_{13}$ & $a^0_{21}$ & $a^0_{32}$ & Exp. & $a^0_{13}$ & $a^0_{21}$ & $a^0_{32}$\\\hline
		1   & 0.50 & 0.50 & 0.50 & 5   & 0.25 & 1.00 & 0.50\\
		2  & 1.00 & 0.50 & 0.25 &  6  & 0.25 & 0.50 & 1.00\\
		3 & 1.00 & 0.25 & 0.50 &   7 & 0.50 & 0.25 & 1.00\\
		4  & 0.50 & 1.00 & 0.25 &     &  &  & \\\hline 
	\end{tabular}
\end{table} 


The results indicate that our conjecture holds, since the best EMP was, in all cases,
the one in which the largest module was a \dm.
Furthermore, as seen in Table \ref{tab:3nodeResults}, in the fully symmetric case
these same EMPs -  VII, IX and XI - outperform the other three - VIII, X and XII.
\begin{table}[!h]
	\centering
	\caption{Trace of the covariance matrix for the cases in Table \ref{tab:3nodeParameters} and
		different EMPs, all equally excited and measured.}
	\label{tab:3nodeResults}
	\begin{adjustbox}{width=\columnwidth}
	\begin{tabular}{c c c c c c c}\hline
		Exp./EMP & VII & VIII & IX & X& XI& XII\\\hline
		1  &\bt0.083 &0.131 &\bt0.083 &0.131 &\bt0.083 &0.131 \\
		2  &0.349 &0.031 &\bt0.027 &0.643 &0.099 &0.101 \\
		3  &0.099 &0.101 &\bt0.027 &0.643 &0.349 &0.031 \\
		4  &0.349 &0.031 &0.099 &0.101 &\bt0.027 &0.643 \\
		5  &0.099 &0.109 &0.349 &0.031 &\bt0.027 &0.643 \\
		6  &\bt0.027 &0.643 &0.349 &0.031 &0.099 &0.109 \\
		7  &\bt0.027 &0.643 &0.099 &0.109 &0.349 &0.031 \\\hline   
	\end{tabular}
	\end{adjustbox}
\end{table}

From these results one concludes that for a three-node network one should
chose an EMP for which the number of inputs is equal to the number of measurements. Specifically, the minimal EMPs VII, IX, XI are more likely to give a best result. Moreover, the best minimal EMPs were approximately 24 times more accurate than the worst minimal EMPs in this Table. 
We have tested numerically these hypotheses for a 5-node cyclic networks and observed that these same principles remain valid.


\subsection{4-node Cyclic Dynamic Networks}
\label{secsec:4nLoop}

Now, consider the case of 4-node cyclic dynamic networks, for which Theorem \ref{teo:CyclicEven} can be applied. In this case, only two minimal EMPs exist, namely I - $\mathcal{B} = \left\lbrace 1, 3\right\rbrace , \mathcal{C} = \left\lbrace 2, 4\right\rbrace $ and II - $\mathcal{B} = \left\lbrace 2, 4\right\rbrace , \mathcal{C} = \left\lbrace 1, 3\right\rbrace $. 

Following the same approach of the previous subsections, we calculate the predictor gradient and analyze the EMPs in order to decide which ones produce the best asymptotic covariance of the parameter estimates. 

In the same spirit of the previous subsections, we analyze several identification scenarios. Considering the fully symmetric case, both EMPs provide the same
overall accuracy, just like in the two-node case, which comes as no surprise.

We conducted a numerical experiment with 10,000 systems where the network's parameters were randomly selected from
a uniform distribution with support in $[-1~, ~1]$. 
The following conjecture was tested:
if $\left[ a^0_{21}a^0_{43}\right]^2 > \left[ a^0_{32}a^0_{14}\right]^2 $ then EMP I is the best, otherwise EMP II.
This conjecture comes from the extrapolation of the rationale used
in the previous networks concerning direct modules: 
in EMP I the modules $a^0_{21}$ and $a^0_{43}$ are the direct modules, so one expects that if they are the larger
modules then this will be the best EMP (mutatis mutandis for EMP II).  
From the $10,000$ systems tested, the conjecture proved correct in 99.64 \% of them, implying
that those indicators can be used to determine which EMP will have the best accuracy for cyclic networks with four nodes.   
Once again it was observed that the gains obtained in choosing the best EMP
in can be very significant: in $20\%$ of the cases the ratio between the variances
of the two EMPs was above $100$, and the median of this ratio was found to be
$8.7$.

Thus, the principle that EMPs where the larger modules are direct modules provide better accuracy is confirmed once again.

\section{Branches}
\label{sec:Branches}

In this section, we deal with dynamic networks with branch structure, which are characterized by a network matrix in the form (\ref{eq:GCasc}).
We will treat the most basic case first, which is a branch with only three nodes, and larger branches are considered later. 

\subsection{Three-nodes}
\label{subsec:3nBranch}

A 3-node dynamic network with branch structure provides only two minimal EMPs, namely:
\begin{enumerate}[(I)]
	\item $\mathcal{B} = \lbrace 1, 2\rbrace$; $\mathcal{C} = \lbrace 3 \rbrace$;   \label{item:B12C3}
	\item 	$\mathcal{B} = \lbrace 1\rbrace$; $\mathcal{C} = \lbrace 2, 3 \rbrace$. \label{item:B1C23}
\end{enumerate}
The difference between the two in whether node 2 is measured or excited. 
Consider first EMP (\ref{item:B12C3}) and
recall that the noise processes obey assumption (\ref{item:noise-lambda}).
The gradient of the optimal predictor can be computed as:
\begin{align}
\psi^{\text{\ref{item:B12C3}}}(t) = 
\begin{bmatrix}
a_{32}r_1(t-2) & a_{21}r_1(t-2) + r_2(t-1)
\end{bmatrix}^T. \nonumber
\end{align}
where the superscript stands for EMP (\ref{item:B12C3}). 
Then, one can obtain the covariance matrix by calculating $P^{\text{\ref{item:B12C3}}} = (\E \psi^{\text{\ref{item:B12C3}}}(t)\psi^{\text{\ref{item:B12C3}}}(t)^T)^{-1}\lambda_3$. 
From this, it is possible to obtain the variance of the parameters' estimates as:
\begin{align}
\text{var}(\hat{a}^{\text{\ref{item:B12C3}}}_{21}) &= \lambda_3([{a_{21}^0]^2}/{\sigma_2^2} + {1}/{\sigma_1^2})/{[a^0_{32}]^2}, \label{eq:varb213node} \\
\text{var}(\hat{a}^{\text{\ref{item:B12C3}}}_{32}) &=	{\lambda_3}/{\sigma_2^2}. 
\end{align} 
Now, consider EMP (\ref{item:B1C23}), in which the second node is excited instead of measured. In this case the gradient of the optimal predictors is:
\begin{align}
\psi^{\text{\ref{item:B1C23}}}(t) = 
\begin{bmatrix}
{r_1(t-1)} & {a_{32}r_1(t-2)} \\
0 & {a_{21}r1(t-2)}
\end{bmatrix}^T, \nonumber
\end{align}
which results in the following expression for the variances of the parameters' estimates:
\begin{align}
\text{var}(\hat{a}^{\text{\ref{item:B1C23}}}_{21}) &= {\lambda_2}/{\sigma_1^2}, \\
\text{var}(\hat{a}^{\text{\ref{item:B1C23}}}_{32}) &=	([{a_{32}^0]^2\lambda_2 + \lambda_3})/([a^0_{21}]^2\sigma_1^2). \label{eq:varb323node} 
\end{align}
Notice that in each EMP there is one direct module, and that the variance
of the direct module does
not depend on the real values of the parameters, but only on the SNR.

Taking the sum of the variances, it is observed that there is a value of $\sigma_2^2$
for which it is is the same for both EMPs: 
\begin{align}
\bar{\sigma}_2^2 = \frac{\lambda_3\sigma_1^2[a^0_{21}]^2[a^0_{32}]^2\left([a^0_{21}]^2 + 1\right) }{\lambda_2[a^0_{21}]^2[a^0_{32}]^2\left( [a^0_{32}]^2 + 1\right) + \lambda_3\left([a^0_{32}]^2 - [a^0_{21}]^2 \right) }. \label{eq:sig2Motivating}
\end{align} 
If $\sigma^2_2 > \bar{\sigma}^2_2$, then EMP (\ref{item:B12C3}) will be more
accurate, otherwise better accuracy will be obtained by EMP (\ref{item:B1C23}). 
This expression gives a clear guideline for exchanging a measurement for an excitation: if node 2 is to be excited, then it must be with a signal of
amplitude larger than $\bar{\sigma}^2_2$, otherwise it is better to measure
node 2 instead of exciting it. 
When both parameters are equal ($a^0_{21} = a^0_{32}$), the expression (\ref{eq:sig2Motivating}) reduces to:
\begin{align}
\sigma_2^2 = ({\lambda_3}/{\lambda_2})\sigma_1^2, \nonumber
\end{align} 
that is, it is better to excite node $2$ (rather than to measure it) if the SNR $\frac{\sigma_2^2}{\lambda_3}$ is larger than
the SNR $\frac{\sigma_1^2}{\lambda_2}$.
%
On the other hand, in an \eem setting the expressions (\ref{eq:varb213node})-(\ref{eq:varb323node}) lead to the following conclusions:
\begin{itemize}
	\item EMP \ref{item:B12C3} results in a smaller trace of the covariance matrix if and only if $[a^0_{21}]^2 < [a^0_{32}]^2$ - that is, the best EMP is the one in which
	the largest module is a direct module;
	\item if $a^0_{21} = a^0_{32}$ then both EMPs results in the same trace of covariance matrix, but $\text{var}(\hat{a}^{\ref{item:B12C3}}_{21}) > \text{var}(\hat{a}^{\ref{item:B1C23}}_{21})$ and $\text{var}(\hat{a}^{\ref{item:B12C3}}_{32}) < \text{var}(\hat{a}^{\ref{item:B1C23}}_{32})$ -  in each EMP the direct module is estimated
	more precisely.
\end{itemize}
These conclusions are also in accordance with the loop results regarding the role played by 
the \dm. 

\subsection{Four-nodes branches}
\label{subsec:4nBranch}

Branches with 4 nodes provide four different minimal EMPs, listed below:
\begin{enumerate}[(I)]
	\item $\mathcal{B} = \left\lbrace 1, 3 \right\rbrace; \mathcal{C} = \left\lbrace 2, 4 \right\rbrace$; \label{item:B13C24}
	\item $\mathcal{B} = \left\lbrace 1, 2\right\rbrace; \mathcal{C} = \left\lbrace 3, 4\right\rbrace$; \label{item:B12C34}
	\item $\mathcal{B} = \left\lbrace 1, 2, 3\right\rbrace; \mathcal{C} = \left\lbrace 4\right\rbrace$; \label{item:B123C4}
	\item $\mathcal{B} = \left\lbrace 1\right\rbrace; \mathcal{C} = \left\lbrace 2, 3, 4\right\rbrace$. \label{item:B1C234}
\end{enumerate} 

Analytical expressions for the variances of the estimates of each parameter in each EMP are given in  Appendix \ref{appendix4} and were used to derive the following results.

Consider first the \eem~case.
If all modules are equal - $a^0_{21} = a^0_{32} = a^0_{43}$ - then EMP II 
yields the smallest variance. In EMP II the number of measures
is the same as the number of excitations, and it is the first two nodes that
are excited, whereas the last two are measured. This principle also applies for larger
branches, as will be seen in the sequel.

In order to verify the role played by the direct modules,
consider the experiments described 
in Table \ref{tab:directModule4nodeCasc}. It is seen in this Table
that those EMPs for which the largest module is not a direct module perform
much worse than the ones where the \dm~is the largest one.

\begin{table}[!h]
	\centering
	\caption{Effect of the direct module in the EMPs in \eem~ setting.}
	\label{tab:directModule4nodeCasc}
	\begin{tabular}{c c c c c c c}
		\hline
		$a^0_{21}$ & $a^0_{32}$ & $a^0_{43}$ & $tr(P^{\text{\ref{item:B13C24}}})$ &  $tr(P^{\text{\ref{item:B12C34}}})$ &  $tr(P^{\text{\ref{item:B123C4}}})$ &  $tr(P^{\text{\ref{item:B1C234}}})$ \\\hline
		20 & 1 & 1 & 0.03 & 2.01 & 4.04 & 0.01 \\
		1  & 1 & 20& 0.03 & 2.01 & 0.01 & 4.04 \\
		1  &20 & 1 & 8.03 & 0.007& 4.02 & 4.02 \\\hline 
	\end{tabular}
\end{table}  

These results are expected, as in the three nodes case, one may find which EMPs are not very accurate by looking to those that do not have their largest modules as direct modules. However, unlike the three node case, there are other factors
to be considered here (and also in larger branches) that influence which EMP
is the best: the cardinalities $|B|$ and $|C|$ and whether the measurements
are made closer to the source of the branch or closer to its sink. In any given 
network these factors will be competing against each other to determine which 
EMP is best. For example, EMP II is best when all modules have the same magnitude,
because it has the best distribution: two measurements at the last two nodes
and two excitations at the first two nodes. But in EMPs I and IV  $a_{21}$ is a direct
module, so in a network where this module is larger than the other ones this
effect may dominate the advantages of EMP II. Indeed this is what is observed
in  the first line of Table \ref{tab:directModule4nodeCasc}.
So, the magnitude of the direct module does not indicate the best EMP by itself,
but it gives  insight into those EMPs that are good candidates, along with the
criteria observed for the fully uniform network.

Let us consider a numeric experiment to see how these different factors
play against each other in determining the best EMP.
We tested a total of 10,000 4-node branch networks with random parameters, where all variables involved (modules and input/noise variances) were selected from a uniform distribution $\mathcal{U}(0, 50)$. In each network the best EMP was selected as the most accurate according to its trace of covariance matrix, resulting in the
distribution shown in Table \ref{tab:EMP2better}.

\begin{table}[!h]
	\caption{How often an EMP was the most accurate for 4-node networks}
	\label{tab:EMP2better}
	\centering
	\begin{tabular}{c c c c c}\hline
		EMP & \ref{item:B13C24}     & \ref{item:B12C34}       & \ref{item:B123C4}       & \ref{item:B1C234}  \\ \hline
		\%  & 1.15\% & 87.94\% & 5.50\% &  5.41\%\\ \hline
	\end{tabular}
\end{table}
These results show that EMP (\ref{item:B12C34}) is the best alternative,
as it was selected as the best in the vast majority of cases.
One should pick this EMP unless some prior on the network suggests otherwise - for
instance, if it is known in advance that $a_{32}$ is significantly
smaller than the other modules.

\subsection{Larger branches}
For branches with more than four nodes the analytical expressions of the variances
can also be obtained in the spirit of the previous sections. However, the large
number of EMPs and of variables to consider (input and noise amplitudes
for each node, modules magnitudes) make it hard to derive direct 
comparisons among the different EMPs directly from the analytical expressions.
So, to check the generalization
of the results obtained so far we have performed
numerical experiments, similar to those described in Table \ref{tab:EMP2better},
for networks of five, six, seven, and eight nodes.
For each cardinality of the network we have
randomly generated networks and selected the best EMP for each network.
The results are shown in Tables \ref{tab:5nbranch} to \ref{tab:8nbranch}. 
The description of these EMPs is given in Appendix \ref{app:empsbranches}.
\begin{table}[!h]
	\caption{How often an EMP was the most accurate for 5-node networks.}
	\label{tab:5nbranch}
	\centering
	\begin{adjustbox}{width=\columnwidth}
	\begin{tabular}{rcrcrcrc}
		\toprule
		EMPs &      \% & EMPs &      \% & EMPs &      \% & EMPs &      \% \\
		\midrule
		1 &   2.02 & 3 &   1.15 & 5 &   0.69 & 7 &   0.68  \\
		2 &  46.10 & 4 &  45.79 & 6 &   1.63 & 8 &   1.94  \\
		\bottomrule
	\end{tabular}
	\end{adjustbox}
\end{table}

\begin{table}[!h]
	\caption{How often an EMP was the most accurate for 6-node networks.}
	\label{tab:6nbranch}
	\centering
	\begin{adjustbox}{width=\columnwidth}
	\begin{tabular}{rcrcrcrc}
		\toprule
		EMPs &     \% & EMPs &     \%  & EMPs &     \% & EMPs &     \% \\
		\midrule
		1 &   0.40 & 5 &   0.27 & 9 &   0.10  & 13 &   0.17   \\
		2 &  16.63 & 6 &   2.39 & 10 &   0.37 & 14 &   0.25   \\
		3 &   0.26 & 7 &   0.34 & 11 &   0.10 & 15 &   0.13   \\
		4 &  60.33 & 8 &  17.66 & 12 &   0.27 & 16 &   0.33   \\
		
		\bottomrule
	\end{tabular}
	\end{adjustbox}
\end{table}

\begin{table}[!h]
	\caption{How often an EMP was the most accurate for 7-node networks.}
	\label{tab:7nbranch}
	\centering
	\begin{adjustbox}{width=\columnwidth}
	\begin{tabular}{rcrcrcrc}
		\toprule
		EMPs &     \% & EMPs &     \%  & EMPs &     \% & EMPs &     \% \\
		\midrule
		1 &   0.08 & 9 &   0.07  & 17 &   0.01 & 25 &   0.01  \\
		2 &   5.53 & 10 &   0.50 & 18 &   0.16 & 26 &   0.03  \\
		3 &   0.03 & 11 &   0.01 & 19 &   0.01 & 27 &   0.02  \\
		4 &  42.19 & 12 &   0.95 & 20 &   0.09 & 28 &   0.06  \\
		5 &   0.01 & 13 &   0.05 & 21 &   0.03 & 29 &   0.01  \\
		6 &   1.34 & 14 &   0.50 & 22 &   0.01 & 30 &   0.05  \\
		7 &   0.10 & 15 &   0.17 & 23 &   0.02 & 31 &   0.01  \\
		8 &  41.96 & 16 &   5.90 & 24 &   0.02 & 32 &   0.07  \\
		\bottomrule
	\end{tabular}
	\end{adjustbox}
\end{table}

\begin{table}[!h]
	\caption{How often an EMP was the most accurate for 8-node networks.}
	\label{tab:8nbranch}
	\centering
	\begin{adjustbox}{width=\columnwidth}
	\begin{tabular}{rcrcrcrc}
		\toprule
		EMPs &     \% & EMPs &     \%& EMPs &     \%& EMPs &     \%\\
		\midrule
		1 &   0.02 & 17 &   0.01 & 33 &   0.00 & 49 &   0.01\\
		2 &   1.69 & 18 &   0.09 & 34 &   0.04 & 50 &   0.00\\
		3 &   0.01 & 19 &   0.00 & 35 &   0.00 & 51 &   0.00\\
		4 &  21.37 & 20 &   0.20 & 36 &   0.02 & 52 &   0.02\\
		5 &   0.01 & 21 &   0.01 & 37 &   0.00 & 53 &   0.01\\
		6 &   0.61 & 22 &   0.04 & 38 &   0.02 & 54 &   0.00\\
		7 &   0.01 & 23 &   0.01 & 39 &   0.01 & 55 &   0.00\\
		8 &  50.09 & 24 &   0.51 & 40 &   0.03 & 56 &   0.03\\
		9 &   0.01 & 25 &   0.00 & 41 &   0.00 & 57 &   0.00\\
		10 &   0.26& 26 &   0.23 & 42 &   0.00 & 58 &   0.01 \\
		11 &   0.01& 27 &   0.02 & 43 &   0.02 & 59 &   0.01 \\
		12 &   1.53& 28 &   0.41 & 44 &   0.01 & 60 &   0.00 \\
		13 &   0.03& 29 &   0.03 & 45 &   0.00 & 61 &   0.01 \\
		14 &   0.33& 30 &   0.15 & 46 &   0.00 & 62 &   0.01 \\
		15 &   0.03& 31 &   0.01 & 47 &   0.01 & 63 &   0.00 \\
		16 &  20.19& 32 &   1.77 & 48 &   0.00 & 64 &   0.04 \\
		
		\bottomrule
	\end{tabular}
	\end{adjustbox}
\end{table}

From the results, both analytical and numerical, obtained for branches of up to four nodes 
one can expect that the best EMP will be the one in which the first half nodes are
excited and the remaining ones are measured. It is seen in Tables \ref{tab:6nbranch} 
and \ref{tab:8nbranch} that this is indeed what happens. 
When this ``half" is not an integer - 
that is, when the number of nodes is odd - it is left do decide
whether to measure or to excite the ``middle" node, and one expects both
choices to be equivalent from the point of view of overall precision. This is what
is observed in Tables \ref{tab:5nbranch} to \ref{tab:7nbranch}, where these two choices are by far the
best, with very little difference between the two.
Moreover, it can be shown analytically that for the \eem scenario 
the variance is the same 
for these two EMPs.

\section{A more general example}
\label{sec:anexample}

\begin{figure}[!h]
	\centering
	\caption{An example of a network with branch and cycle structures combined.}
	\includegraphics[width=0.85\linewidth,keepaspectratio]{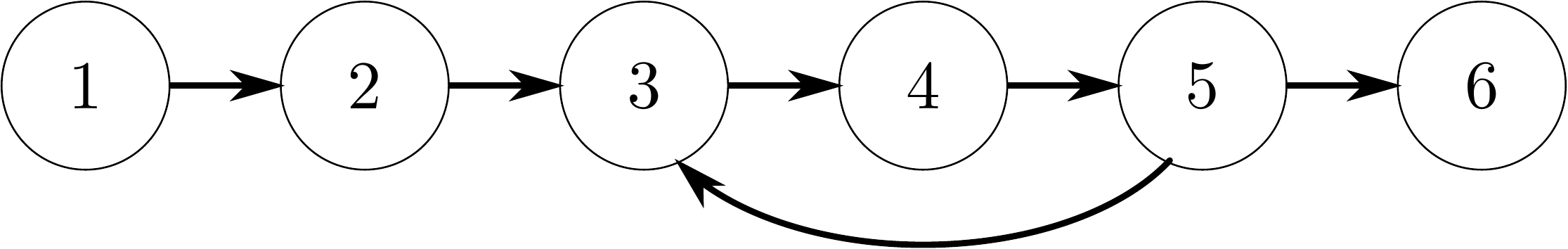}
	\label{fig:complexnet}
\end{figure}

In the previous sections we have separately analyzed branch and cyclic networks. 
Here, we will present an example of a network where these two structures are combined to form a single network. Our objective is to verify to which
extent the results obtained from these two structures can be applied in a more complex topology. 
The network analyzed is depicted in Figure \ref{fig:complexnet}.

From Corollary \ref{cor:summary}, we know that node 1 needs to be excited and node 6 measured. Furthermore, node 2 could be either excited or measured, while the nodes that compose the cycle $\lbrace 3, 4, 5\rbrace$ need to satisfy the conditions of Theorem \ref{teo:IdentLoop}, that is, at least one of them must be excited and measured while the others could be either excited or measured. 
We thus have a total of 24 ``minimal"\footnote{See the remark after Corollary \ref{cor:summary}} EMPs
 which are described in Table \ref{tab:empscomplex} given in Appendix \ref{app:empscomplex}. 

According to the results derived in previous sections, we expect for branches that the first half is excited and the other half is measured. If we consider the cycle (nodes $\{3, 4, 5\}$) as a supernode, it is expected that best accuracy will be obtained for the EMPs whose node 2 is excited. Furthermore, EMPs for which the number of excitations and measurements are equal in the cycle are presumably more accurate; these principles are observed in EMPs 7-12 in Table \ref{tab:empscomplex}. Hence, these EMPs are our candidates to be among the most accurate. On the other hand, the EMPs that most violate these features are 13-18, where node $2$ is measured and the excitation/measurement in the cycle is not equally distributed. For this reason, they are expected to yield poor estimates.  

We have performed a numerical simulation of $1,000$ networks.
The covariance matrices were evaluated numerically. 
As we have done previously, a comparison for all EMPs is first presented for the fully symmetrical network (when every quantity in the nodes is equal). 
A numeric experiment was done with the modules equal to $0.3$, and input/noise variances as $\sigma^2 =1$ and $\lambda = 0.01$, respectively. We have tested all EMPs and the trace of the covariance matrix of each one is given in Table \ref{tab:bcex1}. 
\begin{table}[!h]
	\centering
	\caption{Trace of covariance matrix in the fully symmetrical case for all EMPs.}
	\begin{adjustbox}{width=\columnwidth}
	\begin{tabular}{lrrrrrrr}
		\toprule
		EMP &          trace & EMP &          trace & EMP &          trace & EMP &          trace\\
		\midrule
		1  &   0.755 & 7  &   0.435 &13 &   0.795 & 19 &   0.445 \\
		2  &   1.646 & 8  &   0.631 &14 &   1.448 & 20 &   0.678 \\
		3  &   15.87 & 9  &   0.242 &15 &  13.656 & 21 &   0.290 \\
		4  &   2.484 & 10 &   2.419 &16 &   2.718 & 22 &   2.254 \\
		5  &   0.438 & 11 &   1.515 &17 &   0.958 & 23 &   1.348 \\
		6  &   0.721 & 12 &   0.626 &18 &   1.268 & 24 &   1.019 \\
		\bottomrule
	\end{tabular}
	\end{adjustbox}
	\label{tab:bcex1}
\end{table}

The best EMP (9) is the one for which node two was excited and the number of inputs are equal to the number of outputs in the cycle. This result is expected according to the principles presented for branches and cycles. Furthermore, we notice that the best EMP (9) is approximately 75 times better than the worst EMP (3). Notice that some EMPs (13-18), which violate both principles, are among the worst choices. Interestingly, EMP 21 obeys only the cycle condition and it is the runner-up in this experiment. Although the insights provided in this paper serve as a strong indicator for which EMPs are more accurate, other underlying principles may play a role in the EMPs' precision for more complex networks. For instance, EMP 10 should - according to our principles - yield a good result, while EMP 21 that obeys only the cycle condition is almost as accurate as the best EMP. 

Now, starting from the previous fully symmetrical network, we have tripled the magnitude of one of the modules (resulting in six different networks) and again tested all EMPs in each case to check if the direct module is determinant. The best EMP in each case is given in Table \ref{tab:bcex2}.

\begin{table}[!h]
	\centering
	\caption{Networks with a single module as the largest one}
	\begin{adjustbox}{width=\columnwidth}
	\begin{tabular}{l|rrrrrr}
		\toprule
		largest module & $a_{21}$ & $a_{32}$ & $a_{35}$ & $a_{43}$ & $a_{54}$ & $a_{65}$ \\
		\midrule
		EMP & 21 & 9 & 9 & 8 & 7 & 9 \\
		\midrule
		trace & 0.1607 & 0.1027 & 0.1565  & 0.1383 & 0.1609 & 0.1976\\
		\bottomrule
	\end{tabular}
	\end{adjustbox}
	\label{tab:bcex2}
\end{table} 

In all cases, the most accurate EMP had the largest module as a direct module, which confirms again that this principle should be taken into account when choosing among different EMPs. Remarkably, EMP 21 is the only one where the principles do not hold simultaneously, which indicates that the direct module ($a_{21}$) is decisive for the accuracy of this EMP. 	

Finally, we have tested a thousand random networks with modules drawn from a uniform distribution with range $[-1, 1]$. Table \ref{tab:ntimescomplex} presents the number of times that some EMPs were selected as the best for each simulation. EMPs that are not present there were not selected as the best not even once.

\begin{table}[!h]
	\centering
	\caption{Number of times each EMP was selected as the best  when all modules were randomly selected from $[-1, 1]$.}
	\begin{adjustbox}{width=\columnwidth}
	
	\begin{tabular}{lrrrrrrr}
		\toprule
		EMP &  Freq.  & EMP &    Freq. & EMP &    Freq. & EMP &    Freq.\\
		\midrule
		4  &    9 &9  &  197 & 14 &    2&20 &  125 \\
		5  &   23 &10 &    6 & 16 &    3&21 &  124 \\
		6  &   24 &11 &   13 & 17 &    1&22 &   17 \\
		7  &   27 &12 &  125 & 18 &    1&23 &   36 \\
		8  &  124 &13 &   10 & 19 &   74&24 &   59 \\
		\bottomrule
	\end{tabular}
	\end{adjustbox}
	\label{tab:ntimescomplex}
\end{table}

As can be seen from this Table, three of the candidates EMPs (8, 9 and 12) to have good accuracy were selected in \%44.6 of the simulations. Furthermore, there is a clear advantage for EMP 9, which outperformed the others. The EMPs (13-18) that do not follow neither of the principles were chosen as best EMPs only in few runs. 
We highlight that in 93.7\% of the simulations the largest module was a direct module of the best EMP. This indicates 
once again that direct modules are important factors in the choice of the best EMP.

\section{Conclusion}
\label{sec:Conclusion}

Whereas a lot of literature has been, and is being, produced
on determining which excitation and measurement patterns are viable for
network identification, the question of choosing which of the viable
EMPs is best seems to be fully unexplored;
in this paper we have launched the exploration of this subject. 
We have formalized the problem and proposed a framework for
its study, and we have focused our investigation - an exploratory
one  - to a narrow class of networks
with a specific quality criterion (the trace of the covariance matrix) and a specific
cost to be minimized (the number of nodes involved in the EMP).
Though the results have been obtained for branches and cycles only, 
we have illustrated, by means of an example, how they can be
applied to networks of more general  topology.

From the theoretical analysis and the thousands of randomly generated 
case studies  that have been analyzed, three
principles have emerged for the choice of EMPs. First, it has been noted 
EMPs with equal shares of measurement and excitation - that is, where 
the number of nodes that are excited is equal to the number of nodes that are measured - give the 
best results. Second, in a branch the best results are obtained
by exciting the first nodes (that is, those closer to the source)
and measuring the last ones (those closer to the sink).
And third, the sizes of the  direct modules play a decisive role, so one
should try to choose an EMP such that the larger modules are
direct modules. While the first two principles concern only the topology of
the network, the third one is related to the magnitudes of the modules to be identified.
As such, the third principle requires some prior to be applied, like
knowing which module(s) is(are) of larger magnitude. 
Last, but certainly not least, the huge discrepancy observed
among the different EMPs in the qualities of the estimates 
attests to the paramount relevance of the subject
studied in this paper.

%
%
%
%



\appendix

\section{Covariance Expressions for 2 node cycle}
\label{app1}

The covariance matrices for the 2-node cycle network are given. They can be shown to be as follows:

\begin{align*}
P^{\text{I}} = \frac{\lambda_1}{d_1}
\left[\begin{matrix}
\gamma_0  \left(\sigma_1^{2} + \sigma_2^{2} [a_{12}^0]^2\right)
& -  \left(\gamma_0 \sigma_1^{2} a_{21}^0/a_{12}^0 + \gamma_2 \sigma_2^{2} \right)\\
* &
\gamma_0  \left(\sigma_1^{2} [a_{21}^0]{2} + \sigma_2^{2}\right)/[a^0_{12}]^2
\end{matrix}\right],
\end{align*}

\begin{align*}
P^{\text{II}} = \frac{\lambda_2}{d_2}\left[\begin{matrix}\frac{\gamma_0  \left(\sigma_1^{2} + \sigma_2^{2} [a_{12}^0]^2\right) }{[a_{21}^0]^2 } & - \frac{ \left(\gamma_0 \sigma_2^{2} a_{12}^0 + \gamma_2 \sigma_1^{2} a_{21}^0\right)}{a_{21}^0 }\\ * & \gamma_0  \left(\sigma_1^{2} [a_{21}^0]^2 + \sigma_2^{2}\right)\end{matrix}\right],
\end{align*}

\begin{align*}
P^{\text{III}} = \frac{\lambda_1 \lambda_2}{d_3}\left[\begin{matrix}\frac{\gamma_0  \left(\lambda_1 + \lambda_2 [a_{12}^0]^{2}\right)}{ [a_{21}^0]^{2}} & - \frac{ \left(\gamma_0 \lambda_2 a_{12}^0 + \gamma_2 \lambda_1 a_{21}^0\right)}{a_{21}^0 }\\ * & \gamma_0  \left(\lambda_1 [a_{21}^0]^{2} + \lambda_2\right)\end{matrix}\right],
\end{align*}

\begin{align*}
P^{\text{IV}} = \frac{\lambda_1 \lambda_2}{d_4}\left[\begin{matrix} \gamma_0  \left(\lambda_1 + \lambda_2 [a_{12}^0]^{2}\right) & - \frac{ \left(\gamma_0 \lambda_1 a_{21}^0 + \gamma_2 \lambda_2 a_{12}^0\right)}{ a_{12}^0}\\ * & \frac{\gamma_0  \left(\lambda_1 [a_{21}^0]^{2} + \lambda_2\right)}{[a_{12}^0]^{2}}\end{matrix}\right].
\end{align*}

\section{Variances for the 4-node branch}\label{appendix4}

Here the expressions for modules' variances of a 4-node branch network are given. Applying the same procedure of the branch with three nodes, we can reach the following variance expressions for EMP \ref{item:B13C24}: 
\begin{align}
\text{var}(\hat{a}^{\text{\ref{item:B13C24}}}_{21}) &=  {\lambda_2}/{\sigma_1^{2}},\, \text{var}(\hat{a}^{\text{\ref{item:B13C24}}}_{43}) = {\lambda_4}/{\sigma_3^{2}}, \nonumber \\
\text{var}(\hat{a}^{\text{\ref{item:B13C24}}}_{32}) &=	\frac{\lambda_2 [a^0_{32}]^2}{\sigma_1^{2} [a^0_{21}]^2} + \frac{\lambda_4 [a^0_{32}]^{2}}{\sigma_3^{2} [a^0_{43}]^2} + \frac{\lambda_4}{\sigma_1^{2} [a^0_{21}]^{2} [a^0_{43}]^{2}}. \nonumber 
\end{align}
Similarly, the variance of the parameters in EMP II are as follows: 
\begin{align}
&\text{var}(\hat{a}^{\text{\ref{item:B12C34}}}_{21}) =  {\lambda_3 \lambda_4 \left[[\sigma_1 a^0_{21}]^{2} + \sigma_2^{2}\right]}/{[\sigma_1 \sigma_2 a^0_{32}]^{2} \left[\lambda_3 [a^0_{43}]^{2} + \lambda_4\right]} , \nonumber\\
&\text{var}(\hat{a}^{\text{\ref{item:B12C34}}}_{32}) =	\frac{\lambda_3 \left( \sigma_2^{2} (\lambda_3 [a^0_{43}]^{2} + \lambda_4) + \lambda_4 \sigma_1^{2} [a^0_{21}]^{2}\right)}{\sigma_2^{2} \left[\lambda_3 [a^0_{43}]^{2} [\sigma_1^{2} [a^0_{21}]^{2}  + \sigma_2^{2}] + \lambda_4 [\sigma_1^{2} [a^0_{21}]^{2} + \sigma_2^{2}]\right]}, \nonumber\\
&\text{var}(\hat{a}^{\text{\ref{item:B12C34}}}_{43}) = \left( {\lambda_3 [a^0_{43}]^{2} + \lambda_4}\right) /\left( {[a^0_{32}]^{2} \left(\sigma_1^{2} [a^0_{21}]^{2} + \sigma_2^{2}\right)}\right) .  \nonumber
\end{align}	
Finally, for EMPs III and IV:
\begin{alignat}{2}
\text{var}(\hat{a}^{\text{\ref{item:B123C4}}}_{21}) &=  \frac{\lambda_4 \left(\sigma_1^{2} [a^0_{21}]^{2} + \sigma_2^{2}\right)}{\sigma_1^{2} \sigma_2^{2} [a^0_{32}]^{2} [a^0_{43}]^{2}},
&&\text{var}(\hat{a}^{\text{\ref{item:B1C234}}}_{21}) =  \frac{\lambda_2}{\sigma_1^{2}},\nonumber\\
\text{var}(\hat{a}^{\text{\ref{item:B123C4}}}_{32}) &=	\frac{\lambda_4 [a^0_{32}]^{2}}{\sigma_3^{2} [a^0_{43}]^{2}} + \frac{\lambda_4}{\sigma_2^{2} [a^0_{43}]^{2}}, &&\text{var}(\hat{a}^{\text{\ref{item:B1C234}}}_{32}) =	\frac{\lambda_2 [a^0_{32}]^{2} + \lambda_3}{\sigma_1^{2} [a^0_{21}]^{2}},\nonumber\\
\text{var}(\hat{a}^{\text{\ref{item:B123C4}}}_{43}) &=\frac{\lambda_4}{\sigma_3^{2}},
&&\text{var}(\hat{a}^{\text{\ref{item:B1C234}}}_{43}) = \frac{\lambda_3 [a^0_{43}]^{2} + \lambda_4}{\sigma_1^{2} [a^0_{21}]^{2} [a^0_{32}]^{2}}. \nonumber 
\end{alignat}

\section{EMPs for the more complex example}
\label{app:empscomplex}

The list of the EMPs is given in Table \ref{tab:empscomplex}. 

\begin{table}[!h]
	\centering
	\caption{List of EMPs for the more complex example}
	\begin{adjustbox}{width=\columnwidth}
	\begin{tabular}{lrlr}
		\toprule                                                                         
		N & EMP($\mathcal{B}, \mathcal{C}$) & N & EMP($\mathcal{B}, \mathcal{C}$) \\  
		\midrule

		1  &  $(\lbrace1, 2, 3, 4, 5 \rbrace,   \lbrace3, 6 \rbrace)$     & 13 &  $(\lbrace1, 3, 4, 5 \rbrace$,  $\lbrace2, 3, 6 \rbrace)$  \\
		2  &  $(\lbrace1, 2, 3, 4, 5 \rbrace,   \lbrace4, 6 \rbrace)$     & 14 &  $(\lbrace1, 3, 4, 5 \rbrace$,  $\lbrace2, 4, 6 \rbrace)$  \\
		3  &  $(\lbrace1, 2, 3, 4, 5 \rbrace,   \lbrace5, 6 \rbrace)$     & 15 &  $(\lbrace1, 3, 4, 5 \rbrace$,  $\lbrace2, 5, 6 \rbrace)$  \\
		4  &  $(\lbrace1, 2, 3 \rbrace      ,  \lbrace3, 4, 5, 6 \rbrace)$ & 16 &  $(\lbrace1, 3 \rbrace       ,\lbrace2, 3, 4, 5, 6 \rbrace)$ \\
		5  &  $(\lbrace1, 2, 4 \rbrace      ,  \lbrace3, 4, 5, 6 \rbrace)$ & 17 &  $(\lbrace1, 4 \rbrace       ,\lbrace2, 3, 4, 5, 6 \rbrace)$ \\
		6  &  $(\lbrace1, 2, 5 \rbrace      ,  \lbrace3, 4, 5, 6 \rbrace)$ & 18 &  $(\lbrace1, 5 \rbrace       ,\lbrace2, 3, 4, 5, 6 \rbrace)$ \\
		7  &  $(\lbrace1, 2, 3, 4 \rbrace   ,   \lbrace3, 5, 6 \rbrace)$   & 19 &  $(\lbrace1, 3, 4 \rbrace    ,\lbrace2, 3, 5, 6 \rbrace)$ \\
		8  &  $(\lbrace1, 2, 3, 5 \rbrace   ,   \lbrace3, 4, 6 \rbrace)$   & 20 &  $(\lbrace1, 3, 5 \rbrace    ,\lbrace2, 3, 4, 6 \rbrace)$ \\
		9  &  $(\lbrace1, 2, 4, 5 \rbrace   ,   \lbrace3, 4, 6 \rbrace)$   & 21 &  $(\lbrace1, 4, 5 \rbrace    ,\lbrace2, 3, 4, 6 \rbrace)$ \\
		10 &  $(\lbrace1, 2, 3, 4 \rbrace   ,   \lbrace4, 5, 6 \rbrace)$   & 22 &  $(\lbrace1, 3, 4 \rbrace    ,\lbrace2, 4, 5, 6 \rbrace)$ \\
		11 &  $(\lbrace1, 2, 3, 5 \rbrace   ,   \lbrace4, 5, 6 \rbrace)$   & 23 &  $(\lbrace1, 3, 5 \rbrace    ,\lbrace2, 4, 5, 6 \rbrace)$ \\
		12 &  $(\lbrace1, 2, 4, 5 \rbrace   ,   \lbrace3, 5, 6 \rbrace)$   & 24 &  $(\lbrace1, 4, 5 \rbrace    ,\lbrace2, 3, 5, 6 \rbrace)$ \\
		\bottomrule
	\end{tabular}
	\end{adjustbox}
\label{tab:empscomplex}
\end{table}

\section{Minimal EMPs for branches}
\label{app:empsbranches}

The minimal EMPs for the branch with five to eight nodes are given in tables:

\begin{table}[!h]
	\centering
	\caption{Minimal EMPs for a 5-node branch network.}
	\label{tab:EMPs5branch}
	\begin{tabular}{rlrl}
		\toprule
		N &                                                                          EMP($\mathcal{B}, \mathcal{C}$) & N & EMP($\mathcal{B}, \mathcal{C}$) \\
		\midrule
		1 &  $(\lbrace{1, 2, 3, 4}\rbrace, \lbrace{5}\rbrace)$ & 5 &  $(\lbrace{1, 3, 4}\rbrace, \lbrace{2, 5}\rbrace)$ \\
		2 &  $(\lbrace{1, 2, 3}\rbrace, \lbrace{4, 5}\rbrace)$ & 6 &  $(\lbrace{1, 3}\rbrace, \lbrace{2, 4, 5}\rbrace)$ \\
		3 &  $(\lbrace{1, 2, 4}\rbrace, \lbrace{3, 5}\rbrace)$ & 7 &  $(\lbrace{1, 4}\rbrace, \lbrace{2, 3, 5}\rbrace)$ \\
		4 &  $(\lbrace{1, 2}\rbrace, \lbrace{3, 4, 5}\rbrace)$ & 8 &  $(\lbrace{1}\rbrace, \lbrace{2, 3, 4, 5}\rbrace)$ \\
		\bottomrule
	\end{tabular}
	
\end{table}

\begin{table}[!h]
	\caption{Minimal EMPs for a 6-node branch network.}
	\centering
	\label{tab:EMPs6branch}
	\begin{tabular}{rcrc}
		\toprule
		N &                                                   EMP($\mathcal{B}, \mathcal{C}$) & N &                                                   EMP($\mathcal{B}, \mathcal{C}$) \\
		\midrule
		1 &  $(\lbrace{1, 2, 3, 4, 5}\rbrace, \lbrace{6}\rbrace)$ & 9 &  $(\lbrace{1, 3, 4, 5}\rbrace, \lbrace{2, 6}\rbrace)$  \\
		2 &  $(\lbrace{1, 2, 3, 4}\rbrace, \lbrace{5, 6}\rbrace)$ & 10 &  $(\lbrace{1, 3, 4}\rbrace, \lbrace{2, 5, 6}\rbrace)$ \\
		3 &  $(\lbrace{1, 2, 3, 5}\rbrace, \lbrace{4, 6}\rbrace)$ & 11 &  $(\lbrace{1, 3, 5}\rbrace, \lbrace{2, 4, 6}\rbrace)$ \\
		4 &  $(\lbrace{1, 2, 3}\rbrace, \lbrace{4, 5, 6}\rbrace)$ & 12 &  $(\lbrace{1, 3}\rbrace, \lbrace{2, 4, 5, 6}\rbrace)$ \\
		5 &  $(\lbrace{1, 2, 4, 5}\rbrace, \lbrace{3, 6}\rbrace)$ & 13 &  $(\lbrace{1, 4, 5}\rbrace, \lbrace{2, 3, 6}\rbrace)$ \\
		6 &  $(\lbrace{1, 2, 4}\rbrace, \lbrace{3, 5, 6}\rbrace)$ & 14 &  $(\lbrace{1, 4}\rbrace, \lbrace{2, 3, 5, 6}\rbrace)$ \\
		7 &  $(\lbrace{1, 2, 5}\rbrace, \lbrace{3, 4, 6}\rbrace)$ & 15 &  $(\lbrace{1, 5}\rbrace, \lbrace{2, 3, 4, 6}\rbrace)$ \\
		8 &  $(\lbrace{1, 2}\rbrace, \lbrace{3, 4, 5, 6}\rbrace)$ & 16 &  $(\lbrace{1}\rbrace, \lbrace{2, 3, 4, 5, 6}\rbrace)$ \\
		\bottomrule
	\end{tabular}
\end{table}

\begin{table}[!h]
	\centering
	\caption{Minimal EMPs for a 7-node branch network.}
	\label{tab:EMPs7branch}
	\begin{tabular}{rlrl}
		\toprule
		N &                                                      EMP($\mathcal{B}, \mathcal{C}$) & N &                                                   EMP($\mathcal{B}, \mathcal{C}$) \\
		\midrule
		1 &  $(\lbrace{1, 2, 3, 4, 5, 6}\rbrace, \lbrace{7}\rbrace)$ &17 &  $(\lbrace{1, 3, 4, 5, 6}\rbrace, \lbrace{2, 7}\rbrace)$  \\
		2 &  $(\lbrace{1, 2, 3, 4, 5}\rbrace, \lbrace{6, 7}\rbrace)$ &18 &  $(\lbrace{1, 3, 4, 5}\rbrace, \lbrace{2, 6, 7}\rbrace)$  \\
		3 &  $(\lbrace{1, 2, 3, 4, 6}\rbrace, \lbrace{5, 7}\rbrace)$ &19 &  $(\lbrace{1, 3, 4, 6}\rbrace, \lbrace{2, 5, 7}\rbrace)$  \\
		4 &  $(\lbrace{1, 2, 3, 4}\rbrace, \lbrace{5, 6, 7}\rbrace)$ &20 &  $(\lbrace{1, 3, 4}\rbrace, \lbrace{2, 5, 6, 7}\rbrace)$  \\
		5 &  $(\lbrace{1, 2, 3, 5, 6}\rbrace, \lbrace{4, 7}\rbrace)$ &21 &  $(\lbrace{1, 3, 5, 6}\rbrace, \lbrace{2, 4, 7}\rbrace)$  \\
		6 &  $(\lbrace{1, 2, 3, 5}\rbrace, \lbrace{4, 6, 7}\rbrace)$ &22 &  $(\lbrace{1, 3, 5}\rbrace, \lbrace{2, 4, 6, 7}\rbrace)$  \\
		7 &  $(\lbrace{1, 2, 3, 6}\rbrace, \lbrace{4, 5, 7}\rbrace)$ &23 &  $(\lbrace{1, 3, 6}\rbrace, \lbrace{2, 4, 5, 7}\rbrace)$  \\
		8 &  $(\lbrace{1, 2, 3}\rbrace, \lbrace{4, 5, 6, 7}\rbrace)$ &24 &  $(\lbrace{1, 3}\rbrace, \lbrace{2, 4, 5, 6, 7}\rbrace)$  \\
		9 &  $(\lbrace{1, 2, 4, 5, 6}\rbrace, \lbrace{3, 7}\rbrace)$ &25 &  $(\lbrace{1, 4, 5, 6}\rbrace, \lbrace{2, 3, 7}\rbrace)$  \\
		10 &  $(\lbrace{1, 2, 4, 5}\rbrace, \lbrace{3, 6, 7}\rbrace)$&26 &  $(\lbrace{1, 4, 5}\rbrace, \lbrace{2, 3, 6, 7}\rbrace)$  \\
		11 &  $(\lbrace{1, 2, 4, 6}\rbrace, \lbrace{3, 5, 7}\rbrace)$&27 &  $(\lbrace{1, 4, 6}\rbrace, \lbrace{2, 3, 5, 7}\rbrace)$  \\
		12 &  $(\lbrace{1, 2, 4}\rbrace, \lbrace{3, 5, 6, 7}\rbrace)$&28 &  $(\lbrace{1, 4}\rbrace, \lbrace{2, 3, 5, 6, 7}\rbrace)$  \\
		13 &  $(\lbrace{1, 2, 5, 6}\rbrace, \lbrace{3, 4, 7}\rbrace)$&29 &  $(\lbrace{1, 5, 6}\rbrace, \lbrace{2, 3, 4, 7}\rbrace)$  \\
		14 &  $(\lbrace{1, 2, 5}\rbrace, \lbrace{3, 4, 6, 7}\rbrace)$&30 &  $(\lbrace{1, 5}\rbrace, \lbrace{2, 3, 4, 6, 7}\rbrace)$  \\
		15 &  $(\lbrace{1, 2, 6}\rbrace, \lbrace{3, 4, 5, 7}\rbrace)$&31 &  $(\lbrace{1, 6}\rbrace, \lbrace{2, 3, 4, 5, 7}\rbrace)$  \\
		16 &  $(\lbrace{1, 2}\rbrace, \lbrace{3, 4, 5, 6, 7}\rbrace)$&32 &  $(\lbrace{1}\rbrace, \lbrace{2, 3, 4, 5, 6, 7}\rbrace)$  \\
		\bottomrule
	\end{tabular}
\end{table}

\begin{table}[!h]
	\centering
	\caption{Minimal EMPs for a 8-node branch network.}
	\label{tab:EMPs8branch}
	\begin{adjustbox}{width=\columnwidth}
	\begin{tabular}{rlrl}
		\toprule
		N &                                                          EMP($\mathcal{B}, \mathcal{C}$) & N &                                                   EMP($\mathcal{B}, \mathcal{C}$) \\
		\midrule
		1 &  $(\lbrace{1, 2, 3, 4, 5, 6, 7}\rbrace, \lbrace{8}\rbrace)$ & 33 &  $(\lbrace{1, 3, 4, 5, 6, 7}\rbrace, \lbrace{8, 2}\rbrace)$\\
		2 &  $(\lbrace{1, 2, 3, 4, 5, 6}\rbrace, \lbrace{8, 7}\rbrace)$ & 34 &  $(\lbrace{1, 3, 4, 5, 6}\rbrace, \lbrace{8, 2, 7}\rbrace)$\\
		3 &  $(\lbrace{1, 2, 3, 4, 5, 7}\rbrace, \lbrace{8, 6}\rbrace)$ & 35 &  $(\lbrace{1, 3, 4, 5, 7}\rbrace, \lbrace{8, 2, 6}\rbrace)$\\
		4 &  $(\lbrace{1, 2, 3, 4, 5}\rbrace, \lbrace{8, 6, 7}\rbrace)$ & 36 &  $(\lbrace{1, 3, 4, 5}\rbrace, \lbrace{8, 2, 6, 7}\rbrace)$\\
		5 &  $(\lbrace{1, 2, 3, 4, 6, 7}\rbrace, \lbrace{8, 5}\rbrace)$ & 37 &  $(\lbrace{1, 3, 4, 6, 7}\rbrace, \lbrace{8, 2, 5}\rbrace)$\\
		6 &  $(\lbrace{1, 2, 3, 4, 6}\rbrace, \lbrace{8, 5, 7}\rbrace)$ & 38 &  $(\lbrace{1, 3, 4, 6}\rbrace, \lbrace{8, 2, 5, 7}\rbrace)$\\
		7 &  $(\lbrace{1, 2, 3, 4, 7}\rbrace, \lbrace{8, 5, 6}\rbrace)$ & 39 &  $(\lbrace{1, 3, 4, 7}\rbrace, \lbrace{8, 2, 5, 6}\rbrace)$\\
		8 &  $(\lbrace{1, 2, 3, 4}\rbrace, \lbrace{8, 5, 6, 7}\rbrace)$ & 40 &  $(\lbrace{1, 3, 4}\rbrace, \lbrace{2, 5, 6, 7, 8}\rbrace)$\\
		9 &  $(\lbrace{1, 2, 3, 5, 6, 7}\rbrace, \lbrace{8, 4}\rbrace)$ & 41 &  $(\lbrace{1, 3, 5, 6, 7}\rbrace, \lbrace{8, 2, 4}\rbrace)$\\
		10 &  $(\lbrace{1, 2, 3, 5, 6}\rbrace, \lbrace{8, 4, 7}\rbrace)$& 42 &  $(\lbrace{1, 3, 5, 6}\rbrace, \lbrace{8, 2, 4, 7}\rbrace)$ \\
		11 &  $(\lbrace{1, 2, 3, 5, 7}\rbrace, \lbrace{8, 4, 6}\rbrace)$& 43 &  $(\lbrace{1, 3, 5, 7}\rbrace, \lbrace{8, 2, 4, 6}\rbrace)$ \\
		12 &  $(\lbrace{1, 2, 3, 5}\rbrace, \lbrace{8, 4, 6, 7}\rbrace)$& 44 &  $(\lbrace{1, 3, 5}\rbrace, \lbrace{2, 4, 6, 7, 8}\rbrace)$ \\
		13 &  $(\lbrace{1, 2, 3, 6, 7}\rbrace, \lbrace{8, 4, 5}\rbrace)$& 45 &  $(\lbrace{1, 3, 6, 7}\rbrace, \lbrace{8, 2, 4, 5}\rbrace)$ \\
		14 &  $(\lbrace{1, 2, 3, 6}\rbrace, \lbrace{8, 4, 5, 7}\rbrace)$& 46 &  $(\lbrace{1, 3, 6}\rbrace, \lbrace{2, 4, 5, 7, 8}\rbrace)$ \\
		15 &  $(\lbrace{1, 2, 3, 7}\rbrace, \lbrace{8, 4, 5, 6}\rbrace)$& 47 &  $(\lbrace{1, 3, 7}\rbrace, \lbrace{2, 4, 5, 6, 8}\rbrace)$ \\
		16 &  $(\lbrace{1, 2, 3}\rbrace, \lbrace{4, 5, 6, 7, 8}\rbrace)$& 48 &  $(\lbrace{1, 3}\rbrace, \lbrace{2, 4, 5, 6, 7, 8}\rbrace)$ \\
		17 &  $(\lbrace{1, 2, 4, 5, 6, 7}\rbrace, \lbrace{8, 3}\rbrace)$& 49 &  $(\lbrace{1, 4, 5, 6, 7}\rbrace, \lbrace{8, 2, 3}\rbrace)$ \\
		18 &  $(\lbrace{1, 2, 4, 5, 6}\rbrace, \lbrace{8, 3, 7}\rbrace)$& 50 &  $(\lbrace{1, 4, 5, 6}\rbrace, \lbrace{8, 2, 3, 7}\rbrace)$ \\
		19 &  $(\lbrace{1, 2, 4, 5, 7}\rbrace, \lbrace{8, 3, 6}\rbrace)$& 51 &  $(\lbrace{1, 4, 5, 7}\rbrace, \lbrace{8, 2, 3, 6}\rbrace)$ \\
		20 &  $(\lbrace{1, 2, 4, 5}\rbrace, \lbrace{8, 3, 6, 7}\rbrace)$& 52 &  $(\lbrace{1, 4, 5}\rbrace, \lbrace{2, 3, 6, 7, 8}\rbrace)$ \\
		21 &  $(\lbrace{1, 2, 4, 6, 7}\rbrace, \lbrace{8, 3, 5}\rbrace)$& 53 &  $(\lbrace{1, 4, 6, 7}\rbrace, \lbrace{8, 2, 3, 5}\rbrace)$ \\
		22 &  $(\lbrace{1, 2, 4, 6}\rbrace, \lbrace{8, 3, 5, 7}\rbrace)$& 54 &  $(\lbrace{1, 4, 6}\rbrace, \lbrace{2, 3, 5, 7, 8}\rbrace)$ \\
		23 &  $(\lbrace{1, 2, 4, 7}\rbrace, \lbrace{8, 3, 5, 6}\rbrace)$& 55 &  $(\lbrace{1, 4, 7}\rbrace, \lbrace{2, 3, 5, 6, 8}\rbrace)$ \\
		24 &  $(\lbrace{1, 2, 4}\rbrace, \lbrace{3, 5, 6, 7, 8}\rbrace)$& 56 &  $(\lbrace{1, 4}\rbrace, \lbrace{2, 3, 5, 6, 7, 8}\rbrace)$ \\
		25 &  $(\lbrace{1, 2, 5, 6, 7}\rbrace, \lbrace{8, 3, 4}\rbrace)$& 57 &  $(\lbrace{1, 5, 6, 7}\rbrace, \lbrace{8, 2, 3, 4}\rbrace)$ \\
		26 &  $(\lbrace{1, 2, 5, 6}\rbrace, \lbrace{8, 3, 4, 7}\rbrace)$& 58 &  $(\lbrace{1, 5, 6}\rbrace, \lbrace{2, 3, 4, 7, 8}\rbrace)$ \\
		27 &  $(\lbrace{1, 2, 5, 7}\rbrace, \lbrace{8, 3, 4, 6}\rbrace)$& 59 &  $(\lbrace{1, 5, 7}\rbrace, \lbrace{2, 3, 4, 6, 8}\rbrace)$ \\
		28 &  $(\lbrace{1, 2, 5}\rbrace, \lbrace{3, 4, 6, 7, 8}\rbrace)$& 60 &  $(\lbrace{1, 5}\rbrace, \lbrace{2, 3, 4, 6, 7, 8}\rbrace)$ \\
		29 &  $(\lbrace{1, 2, 6, 7}\rbrace, \lbrace{8, 3, 4, 5}\rbrace)$& 61 &  $(\lbrace{1, 6, 7}\rbrace, \lbrace{2, 3, 4, 5, 8}\rbrace)$ \\
		30 &  $(\lbrace{1, 2, 6}\rbrace, \lbrace{3, 4, 5, 7, 8}\rbrace)$& 62 &  $(\lbrace{1, 6}\rbrace, \lbrace{2, 3, 4, 5, 7, 8}\rbrace)$ \\
		31 &  $(\lbrace{1, 2, 7}\rbrace, \lbrace{3, 4, 5, 6, 8}\rbrace)$& 63 &  $(\lbrace{1, 7}\rbrace, \lbrace{2, 3, 4, 5, 6, 8}\rbrace)$ \\
		32 &  $(\lbrace{1, 2}\rbrace, \lbrace{3, 4, 5, 6, 7, 8}\rbrace)$& 64 &  $(\lbrace{1}\rbrace, \lbrace{2, 3, 4, 5, 6, 7, 8}\rbrace)$ \\
		\bottomrule
	\end{tabular}
	\end{adjustbox}
\end{table}


\bibliographystyle{plain}        
\bibliography{dynetoptbib}           

\end{document}